\documentclass[11pt]{amsart}
\usepackage[all]{xy}
\usepackage{ifthen}
\usepackage{soul}
\usepackage{latexsym}
\usepackage{amsthm,amssymb}
\usepackage{amsbsy,amsfonts,amsmath}
\usepackage{enumitem}
\usepackage{a4}
\usepackage[ansinew]{inputenc}
\usepackage{graphicx,fancyhdr}

\usepackage{multicol}
\usepackage{multienum}

\usepackage[dvips, dvipsnames, usenames]{color}
\allowdisplaybreaks

\newtheorem{theorem}{Theorem}[section]

\newtheorem{lemma}[theorem]{Lemma}

\newtheorem{prop}[theorem]{Proposition}

\theoremstyle{definition}
\newtheorem{definition}[theorem]{Definition}
\newtheorem{example}[theorem]{Example}

\theoremstyle{remark}
\newtheorem{remark}[theorem]{Remark}

\newcommand{\at}[1]{^{(#1)}}
\newcommand{\atb}[1]{^{[#1]}}
\newcommand{\Di}[1]{\operatorname{Dist}{#1}}
\newcommand{\Dis}[2]{\operatorname{Dist}_{#1}{#2}}
\newcommand{\Dl}[2]{\mathcal{D}_{\lambda,{#1}}({#2})}
\newcommand{\Dp}[2]{\mathcal{D}_{\lambda,{#1}}^{+}({#2})}
\newcommand{\Dm}[2]{\mathcal{D}_{\lambda,{#1}}^{-}({#2})}
\newcommand{\Dpm}[2]{\mathcal{D}_{\lambda,{#1}}^{\pm}({#2})}
\newcommand{\Dz}[2]{\mathcal{D}_{\lambda,{#1}}^{0}({#2})}
\newcommand{\Dge}[2]{\mathcal{D}_{\lambda,{#1}}^{\ge0}({#2})}
\newcommand{\Dle}[2]{\mathcal{D}_{\lambda,{#1}}^{\le0}({#2})}

\newcommand{\Dli}[1]{\mathcal{D}_{\lambda}({#1})}

\newcommand{\verma}[2]{\mathcal{M}_{{#1}}({#2})}
\newcommand{\irr}[2]{\mathcal{L}_{{#1}}({#2})}

\newcommand\Ur{\operatorname{U}^{[p]}}

\newcommand\C{\mathbb C}

\newcommand\Z{\mathbb Z}
\newcommand\N{\mathbb N}

\newcommand\g{\mathfrak{g}}

\newcommand\Sl{\mathfrak{sl}}

\def\F{\mathfrak{F}}

\def\Et{\mathtt{E}}
\def\Ft{\mathtt{F}}
\def\Kt{\mathtt{K}}
\def\rt{\mathtt{r}}

\def\cF{\mathcal{F}}

\def\cL{\mathcal{L}}
\def\cN{\mathcal{N}}
\def\cM{\mathcal{M}}

\def\cS{\mathcal{S}}
\def\cD{\mathcal{D}}

\newcommand\ad{\operatorname{ad}}
\newcommand\Ad{\operatorname{Ad}}
\newcommand{\eps}{\varepsilon}

\newcommand\ch{\operatorname{char}}

\newcommand\End{\operatorname{End}}

\newcommand{\Ind}[2]{\operatorname{Ind}_{#1}^{#2}}
\newcommand\SL{\operatorname{SL}}
\newcommand\id{\operatorname{id}}

\newcommand\co{\operatorname{co}}

\newcommand\ind{\operatorname{ind}}
\newcommand\soc{\operatorname{soc}}

\newcommand{\bil}[2]{\begin{Bmatrix} #1 \\ #2 \end{Bmatrix}_{\lambda}}
\newcommand{\bik}[2]{\begin{Bmatrix} #1 \\ #2 \end{Bmatrix}}
\newcommand{\bic}[2]{\begin{bmatrix} #1 \\ #2 \end{bmatrix}_{\lambda}}
\newcommand{\fac}[1]{\begin{bmatrix} #1 \end{bmatrix}_{\lambda}^{!}}
\newcommand{\bra}[1]{\begin{bmatrix} #1 \end{bmatrix}_{\lambda}}
\def\ot{\otimes}
\def\u{\mathfrak{u}}

\def\pf{\begin{proof}}
\def\epf{\end{proof}}

\begin{document}

\title[A quantum version of the algebra of distributions of $\SL_2$]{A quantum version of the algebra of distributions of $\SL_2$}
\author[Iv\'an Angiono]{Iv\'an Angiono}
\address{FaMAF-CIEM (CONICET), Universidad Nacional de C\'ordoba.}
\email{angiono@famaf.unc.edu.ar}

\thanks{ Work partially supported by CONICET,
Secyt (UNC), the MathAmSud project GR2HOPF and ANPCyT (Foncyt), partially done during a visit to the MPI (Bonn), supported by the Alexander von Humboldt Foundation.}

\begin{abstract}
Let $\lambda$ be a primitive root of unity of order $\ell$. We introduce a family of finite-dimensional algebras
$\{\Dl{N}{\Sl_2}\}_{N\in\N_0}$ over the complex numbers, such that $\Dl{N}{\Sl_2}$ is a subalgebra of $\Dl{M}{\Sl_2}$ if $N<M$,
and $\Dl{N-1}{\Sl_2}\subset \Dl{N}{\Sl_2}$ is a $\u_{\lambda}(\Sl_2)$-cleft extension.

The simple $\Dl{N}{\Sl_2}$-modules $(\irr{N}{p})_{0\le p<\ell^{N+1}}$ are highest weight modules, which admit a tensor product
decomposition: the first factor is a simple $\u_{\lambda}(\Sl_2)$-module and the second factor is a simple $\Dl{N-1}{\Sl_2}$-module.
This factorization resembles the corresponding Steinberg decomposition, and the family of algebras resembles the presentation of algebra of distributions of
$\SL_2$ as a filtration by finite-dimensional subalgebras.
\end{abstract}

\maketitle

\section{Introduction}

A difficult question regarding the simple modules over a simple, simply connected algebraic group $G$ over an algebraically closed field of positive characteristic $\Bbbk$ is to find an explicit formula for their characters. A formula involving the action of the corresponding affine Weyl group was proposed by Lusztig \cite{L - problems} in 1980. Subsequently this formula was shown to hold in large characteristic by the combined efforts of Kazhdan-Lusztig, Kashiwara-Tanisaki, Lusztig and Andersen-Jantzen-Soergel. More recently Williamson \cite{W - counterexample} found many counter-examples to the expected bounds in this conjecture. A different approach to character formula with emphasis in the Steinberg decomposition for algebraic groups is given in \cite{L - character}.

Around 1990 Lusztig started to study quantum groups $U_\lambda(\g)$ at a primitive root of unity $\lambda$ of order $\ell$ in order to have algebras over the complex numbers whose representation theory resembles those of simply connected semisimple algebraic groups over algebraically
closed fields of positive characteristic. In particular he conjectured a similar formula for the character of simple modules \cite{L - mod rep}, which holds in this case by a hard proof of Kazhdan-Lusztig. A remarkable fact about $U_\lambda(\g)$ is that it fits into a Hopf algebra extension of the corresponding small quantum group $\u_\lambda(\g)$ by the enveloping algebra $U(\g)$; each simple module satisfies a kind of Steinberg decomposition: it is written as a tensor product of a simple module of $\u_\lambda(\g)$ with a simple module $U(\g)$, viewed as $U_\lambda(\g)$-module via a (kind of) Frobenius map.

A fundamental difference, however, between the representation theory of the algebraic group and the corresponding quantum group at a root of unity is the form of the Steinberg (resp. Lusztig) tensor product theorem: for the algebraic group the theorem involves an arbitrary number of iterations of the Frobenius twist, whereas for the quantum group only one Frobenius twist occurs. It has been proposed by Soergel and Lusztig that there might exist analogues of the quantum group which parallel to a greater and greater extent the representation theory of the algebraic group. Such an object has the potential to deepen our understanding of the representation theory of algebraic groups \cite{W - notes}.

The purpose of this paper is to propose such an object for $\Sl_2$. More precisely, we introduce a family of finite dimensional algebras $\{\Dl{N}{\Sl_2}\}_{N\in\N_0}$ over the complex numbers
to mimic the filtration of the algebra of distributions of $\SL_2$ as a filtration by finite-dimensional subalgebras. This filtration is deeply motivated by the approach proposed in \cite{L - character}. The main objective is to
find a $\C$-algebra whose representation category \emph{behaves} as that of simple, simply connected algebraic groups over algebraically
closed fields of positive characteristic, even more similar than $U_\lambda(\g)$.

\begin{itemize}[leftmargin=*]\renewcommand{\labelitemi}{$\diamond$}
  \item Each algebra $\Dl{N}{\Sl_2}$ is presented by generators and relations; the relations in Definition \ref{defn:Dlambda} resemble those defining finite dimensional subalgebras of the algebra of distributions of $\SL_2$ \cite{Tak - hyperalgebras}.
  \item The first step corresponds to the small quantum group: $\Dl{0}{\Sl_2}\simeq \u_\lambda(\Sl_2)$. If $M<N$, then $\Dl{M}{\Sl_2}$ is a subalgebra of $\Dl{N}{\Sl_2}$, and at the same time there exists a surjective map $\Dl{N}{\Sl_2} \twoheadrightarrow \Dl{M}{\Sl_2}$, see Lemma \ref{lemma:proj over uqsl2}. Thus there exists a surjective map  $\Dl{N}{\Sl_2} \twoheadrightarrow \u_\lambda(\Sl_2)$, a kind of \emph{Frobenius map}.
  \item For the algebra of distributions, there exist extensions of Hopf algebras between consecutive terms of a filtration by (finite dimensional) Hopf subalgebras, see Proposition \ref{prop:extn Hopf algs}. In this case, $\Dl{N-1}{\Sl_2} \subset \Dl{N}{\Sl_2}$ is a
  $\u_\lambda(\Sl_2)$-cleft extension for all $N\in\N$, see Proposition \ref{prop:cleft-ext}.
  \item Each $\Dl{N}{\Sl_2}$ admits a \emph{triangular decomposition} into a positive, a zero and a negative part, see Proposition \ref{prop:basis, triang decomp}. Reasonably each simple module is a \emph{highest weight module}, see Proposition \ref{prop:simple-modules}.
  \item Each simple $\Dl{N}{\Sl_2}$-module admits a \emph{Steinberg decomposition} as the tensor product of a simple $\Dl{N-1}{\Sl_2}$-module and a simple $\u_\lambda(\Sl_2)$-module as stated in Theorem \ref{thm:simple-modules-tensor-decomp}.
\end{itemize}
The next step is to define families of algebras $\Dl{N}{\g}$ for any semisimple Lie algebra $\g$ \cite{An}. This will be the content of a forthcoming paper, where the first step includes the definition of \emph{Lusztig's isomorphisms at level $N$}: that is, to consider Lusztig's isomorphisms for small quantum groups, which induce maps for $\Dl{N}{\g}$. Hence $\Sl_2$ is a key step to prove the existence of a PBW basis labeled by the roots of $\g$.
As the simple modules of these algebras satisfy a Steinberg tensor product decomposition, we hope to attack the modular case from the approach established by Lusztig in \cite{L - character}.

\subsection{Notation}
Let $H$ be a Hopf algebra with counit $\epsilon$ and antipode $\cS$. $H^+$ is the augmentation ideal, i.-e. the kernel of $\epsilon$.
The left adjoint action of $H$ on itself is
$\Ad(a)b=a_1 b \cS(a_2)$, $a,b\in H$. A Hopf subalgebra $A$ is (left) normal if it is stable by
the (left) adjoint action. Given $\pi:H\to K$ a Hopf algebra map,
\begin{align*}
H^{\co \pi} &= \{h\in H: (\id \ot \pi)\Delta(h)=h\ot 1\}, \
^{\co \pi} H = \{h\in H: (\pi \ot \id)\Delta(h)=1\ot h\}, 
\end{align*}
are the sets of left, respectively right, coinvariant elements.

\medskip

Let $A$ be a right $H$-comodule algebra; that is, an $H$-comodule such that the coaction map $\rho:A\to A\ot H$ is an algebra map. Let
\begin{align*}
B:=A^{\co H} &= \{a\in A: \rho(h)=a\ot 1\}, 
\end{align*}
the subalgebra of coinvariant elements. Then $B\subset A$
is a cleft extension if there exists an $H$-colinear convolution-invertible map
$\gamma:H\to A$; we refer to \cite[Section
7]{Mo} for more information.

A sequence of Hopf algebra maps
\begin{align*}
\xymatrix@C=9.5pt{\Bbbk \ar@{->}[rr] & & A \ar@{->}[rr]^{\iota} & & C \ar@{->}[rr]^{\pi} & & B \ar@{->}[rr] & & \Bbbk  }
\end{align*}
is \emph{exact} \cite{And - ext,ADev} if the following conditions hold:
\begin{align*}
\iota&\text{ is injective,} & \pi&\text{ is surjective,}
& \ker\pi &= C \iota(A^{+}), & \iota(A) &= C^{\co \pi}.
\end{align*}

\section{Algebras of Distributions of reductive groups}\label{section:dist alg}

Let $\Bbbk$ be an algebraically closed field, $p=\ch \Bbbk$. Let $G$ be a simply connected semisimple algebraic group with Cartan matrix $A=(a_{ij})_{1\le i,j \le\theta}$,
$\g=\operatorname{Lie} G$. 
Following \cite[Chapter 7]{Jan - book} we recall some definitions concerning algebraic groups. Then we give some results which illustrate those results we want to mimic for the quantum counterpart studied later.

\subsection{The algebra of distributions}
Let $I_e=\{f\in\Bbbk[G]| f(e)=0\}$. A \emph{distribution} on $G$ (with support on $e$) of order $n\in\N_0$ is a linear map $\mu:\Bbbk[G]\to\Bbbk$ such that 
$\mu_{|I_e^{n+1}}\equiv 0$. Let $\Dis{n}{G}$ be the set of all distributions of order $n$, which is a $\Bbbk$-vector space. Now
\begin{align*}
\Di{G}= \bigcup_{n\geq 0} \Dis{n}{G} = \{ \mu \in \Bbbk[G]^{\ast} | \mu_{|I_e^{n+1}}\equiv 0 \text{ for some } n\in\N\}
\end{align*}
is the set of all distributions, which is a $\Bbbk[G]$-module. Then $\Di{G}$ is a Hopf subalgebra of $\Bbbk[G]^{\ast}$, called the \emph{algebra of distributions} (or the \emph{hyperalgebra}) of $G$. As algebra and coalgebra, it is filtered:
\begin{align*}
\Dis{m}{G} \cdot \Dis{n}{G} &\subset \Dis{m+n}{G} & 
\text{for all } & m,n \in\N_0.
\end{align*}

We describe now two basic examples. We refer to \cite[7.8]{Jan - book} for more details.

\begin{example}\label{ex:additive-group}
Let $G=G_a$ (the additive group, which is the spectrum of $\Bbbk[t]$). For each $n\in\N_0$ we set $\gamma_n\in \Bbbk[G_a]^{\ast}$ as the function such that $\gamma_n(t^m)=\delta_{n,m}$ for all $m\in\N_0$. Then $(\gamma_n)_{n\in\N_0}$ is a basis of $\Di{G_a}$, and
\begin{align*}
\gamma_m \gamma_n & = \binom{m+n}{m} \gamma_{m+n}, & \text{for all } & m,n \in\N_0.
\end{align*}
In particular, $\gamma_1^p =p! \gamma_p=0$.
\end{example}

\begin{example}\label{ex:multiplicative-group}
Let $G=G_m$ (the multiplicative group, which is the spectrum of $\Bbbk[t,t^{-1}]$). For each $n\in\N_0$ we set $\varpi_n\in \Bbbk[G_a]^{\ast}$ as the function such that $\varpi_n((t-1)^m)=\delta_{n,m}$ for all $m\leq n$, $\varpi_n(I_e^{n+1})=0$. Then $(\varpi_n)_{n\in\N_0}$ is a basis of $\Di{G_m}$, and the multiplication satisfies that
\begin{align*}
\varpi_m \varpi_n & = \sum_{i=0}^{\min\{m,n\}} \frac{(m+n-i)!}{(m-i)!(n-i!i!)} \varpi_{m+n-1}, & \text{for all } & m,n \in\N_0.
\end{align*}
\end{example}

\medskip

By \cite{Tak - hyperalgebras} the algebra $\Di{G}$ is presented by generators $H_i\at{n}$, $X_i\at{n}$, $Y_i\at{n}$, $1\le i\le\theta$,
$n\in\N_0$, where $H_i\at{0} = X_i\at{0}=Y_i\at{0}=1$, and relations
\begin{align}\label{eq:rels dist-alg -1}
H_i(t) H_i(u) &= H_i(t+u+tu), \\ \label{eq:rels dist-alg -2}
H_i(t) H_j(u) &= H_j(u)H_i(t), \\ \label{eq:rels dist-alg -3}
X_i(t) X_i(u) &= X_i(t+u), \\ \label{eq:rels dist-alg -4}
Y_i(t) Y_i(u) &= Y_i(t+u), \\ \label{eq:rels dist-alg -5}
X_i(t) Y_i(u) &= Y_i\left(\frac{u}{1+tu}\right) H_i(tu) X_i\left(\frac{t}{1+tu}\right), \\ \label{eq:rels dist-alg -6}
X_i(t) Y_j(u) &= Y_j(u)X_i(t), \\ \label{eq:rels dist-alg -7}
H_i(t) X_j(u) &= X_j\left( (1+t)^{a_{ij}} u \right) H_i(t), \\ \label{eq:rels dist-alg -8}
H_i(t) Y_j(u) &= Y_j\left( (1+t)^{-a_{ij}} u \right) H_i(t), \\ \label{eq:rels dist-alg -9}
\ad \left( X_i\at{n}\right)\left( X_j\at{m}\right) &= \sum_{k=0}^n (-1)^k X_i\at{n-k}X_j\at{m}X_i\at{k}=0, & n &>-ma_{ij}, \\ \label{eq:rels dist-alg -10}
\ad \left( Y_i\at{n}\right)\left( Y_j\at{m}\right) &= \sum_{k=0}^n (-1)^k Y_i\at{n-k}Y_j\at{m}Y_i\at{k}=0, & n &>-ma_{ij}.
\end{align}
for $1\le i\neq j \le \theta$, where we consider the following elements of $\Di{G}[[t]]$:
\begin{align*}
 H_i(t) &= \sum_{n=0}^\infty t^n H_i\at{n}, & X_i(t) &= \sum_{n=0}^\infty t^n X_i\at{n}, & Y_i(t) &= \sum_{n=0}^\infty t^n Y_i\at{n}.
\end{align*}

From \eqref{eq:rels dist-alg -2} we have $H_i\at{m}H_j\at{n} = H_j\at{n}H_i\at{m}$ for $i\neq j$, and from \eqref{eq:rels dist-alg -1},
\begin{align}\label{eq:reln Hn Hm}
H_i\at{m}H_i\at{n} &= \sum\limits_{\ell=0}^{\min\{m,n\} } \binom{m+n-\ell}{m}\binom{m}{\ell} H_i\at{m+n-\ell}.
\end{align}
From \eqref{eq:rels dist-alg -3} and \eqref{eq:rels dist-alg -4},
\begin{align}\label{eq:reln Xn Xm}
X_i\at{m}X_i\at{n} &= \binom{m+n}{m} X_i\at{m+n}, & Y_i\at{m}Y_i\at{n} &= \binom{m+n}{m} Y_i\at{m+n}.
\end{align}
From these formulas, the $H_i\at{n}$'s generate a copy of $\Di{G_m}$, while the $X_i\at{n}$'s, respectively the $Y_i\at{n}$'s, generate a copy of $\Di{G_a}$, see Examples \ref{ex:additive-group} and \ref{ex:multiplicative-group}.

From \eqref{eq:rels dist-alg -6} we have $X_i\at{m}Y_j\at{n} = Y_j\at{n}X_i\at{m}$ for $i\neq j$, and from \eqref{eq:rels dist-alg -5},
\begin{align*}
 \sum_{n,m} t^n u^m X_i\at{n}Y_i\at{m} &= \sum_{a,b,c} \frac{u^{a+b}t^{b+c}}{(1+tu)^{a+c}}  Y_i\at{a}H_i\at{b}X_i\at{c} \\
 &= \sum_{a,b,c,d} (-1)^d \binom{a+c+d}{d} u^{a+b+d}t^{b+c+d}  Y_i\at{a}H_i\at{b}X_i\at{c}.
\end{align*}
Thus,
\begin{align}\label{eq:reln Xn Ym}
 X_i\at{n}Y_i\at{m} = \sum_{\ell=0}^{\min\{m,n\} } \sum_{k=0}^\ell (-1)^{\ell-k} \binom{m+n-\ell-k}{\ell-k} Y_i\at{m-\ell}H_i\at{k}X_i\at{n-\ell}.
\end{align}
A particular case of this formula is the following
\begin{align}\label{eq:reln Xn Ym, case}
  [X_i\at{p^n},Y_i\at{p^m}] &= \sum_{\ell=1}^{\min \{p^m, p^n\}}  Y_i\at{p^m-\ell} \left( \sum_{k=0}^\ell \binom{\ell+k}{\ell-k} H_i\at{k} \right) X_i\at{p^n-\ell}.
\end{align}
From \eqref{eq:rels dist-alg -7} and \eqref{eq:rels dist-alg -8} we have
\begin{align}\label{eq:reln Hn Xm}
[H_i\at{p^m}, X_j\at{p^n}] & =  \delta_{n,m} a_{ij} X_i\at{p^n}, &
[H_i\at{p^m}, Y_j\at{p^n}] & =  -\delta_{n,m} a_{ij} Y_i\at{p^n}.
\end{align}

\medskip

\subsection{Hopf algebra extensions from $\Di{G}$}
Let $\cD_n G:=\Dis{p^n}G$. As a consequence of these formulas we have the following result.

\begin{lemma}\label{lemma:normal Hopf subalg}
For all $n\in\N$, $\cD_{n}{G}$ is a normal Hopf subalgebra of $\cD_{n+1}{G}$.
\end{lemma}
\pf
$\cD_{n+1}{G}$ is generated as an algebra by $X_i\at{p^k}$, $Y_i\at{p^k}$, $H_i\at{p^k}$, $0\le k\le n$, so it is enough to prove
that $\cD_{n}{G}$ is stable by the adjoint action of $X_i\at{p^n}$, $Y_i\at{p^n}$, $H_i\at{p^n}$ since the remaining generators
belong to $\cD_{n}{G}$, and $\cD_{n}{G}$ is a Hopf subalgebra.

As $X_i\at{p^n}$ is primitive, $\Ad X_i\at{p^n} = \ad X_i\at{p^n}$. If $m<n$,
then $\ad (X_i\at{p^n})Y_i\at{p^m}$, $\ad (X_i\at{p^n})H_i\at{p^m}\in \cD_{n}{G}$, and $\ad (X_i\at{p^n})X_i\at{p^m}=0$. Let $j\neq i$. Note that
$\ad (X_i\at{p^n})X_j\at{p^m}=\ad (X_i\at{p^n})Y_j\at{p^m}=0$ since they commute, and from \eqref{eq:reln Hn Xm}, $\ad (X_i\at{p^n})H_j\at{p^m}=0$.
Therefore $\ad (X_i\at{p^n})\cD_{n}{G}\subset \cD_{n}{G}$
Analogous computations show that $\ad (Y_i\at{p^n})\cD_{n}{G}$, $\ad (H_i\at{p^n})\cD_{n}{G}\subset \cD_{n}{G}$.
\epf

Now define $\pi_k:\cD_{k+1}{G}\to \cD_{1}{G}=\Ur(\g)$ as follows
\begin{align}\label{eq:defn pi n}
\pi_k(X_i\at{n}) &= \begin{cases} X_i\at{n'}, & \mbox{if }n=p^k n', \\ 0, & \mbox{otherwise}, \end{cases} &
\pi_k(H_i\at{n}) &= \begin{cases} H_i\at{n'}, & \mbox{if }n=p^k n', \\ 0, & \mbox{otherwise}. \end{cases} \\ \notag
\pi_k(Y_i\at{n}) &= \begin{cases} Y_i\at{n'}, & \mbox{if }n=p^k n', \\ 0, & \mbox{otherwise}. \end{cases}
\end{align}

\begin{remark}\label{rem:comb numbers non zero}
Let $n<p^{k+1}$, $0<t<p$. If $p^k$ does not divide $n$, then $\binom{p^k t}{n}=0$, otherwise $n=p^k n'$ and $\binom{p^k t}{n}=\binom{t}{n'}$
by Lucas' Theorem.
\end{remark}

\begin{lemma}\label{lemma:surj Hopf alg map}
$\pi_k$ is a surjective Hopf algebra map.
\end{lemma}
\pf
First we have to check that $\pi_k$ is well defined; i.-e. that the map defined from the free algebra on generators $X_i\at{n}$,
$Y_i\at{n}$ and $H_i\at{n}$ annihilates the defining relations. We check easily that for all $i\neq j$, and $m,n\in \N_0$,
\begin{align*}
\pi_k(H_i\at{m}H_j\at{n} -H_j\at{n}H_i\at{m})&=\pi_k(X_i\at{m}Y_j\at{n} - Y_j\at{n}X_i\at{m})=0.
\end{align*}
For \eqref{eq:reln Hn Hm}, $\pi_k$ annihilates both sides of the equation if $p^k$ does not divides $m$ since either
$\pi_k(H_i\at{m+n-\ell})=0$ or else $\binom{m+n-\ell}{m}=0$. Now set $m=p^k m'$, $n=p^k n'$
\begin{align*}
H_i\at{m}H_i\at{n} &= \sum\limits_{\ell=0}^{\min\{m,n\} } \binom{m+n-\ell}{m}\binom{m}{\ell} H_i\at{m+n-\ell} \\
&= \sum\limits_{\ell'=0}^{\min\{m',n'\} } \binom{p^k(m'+n'-\ell')}{p^km'}\binom{p^k m'}{p^k\ell'} H_i\at{p^k(m'+n'-\ell')} \\
&= \sum\limits_{\ell'=0}^{\min\{m',n'\} } \binom{m'+n'-\ell'}{m'}\binom{m'}{\ell'} H_i\at{p^k(m'+n'-\ell')},
\end{align*}
since $\binom{m}{\ell}=0$ when $p^k$ does not divide $\ell$, so $\pi_k$ applies \eqref{eq:reln Hn Hm} to 0.

For \eqref{eq:reln Xn Xm}, if $p^k$ does not divide $m+n$, then both sides of the equality are annihilated by $\pi_k$.
If $p^k$ divides $m+n$ but does not divide $m$, then again $\pi_k$ annihilates both sides
of \eqref{eq:reln Xn Xm} since $\binom{m+n}{m}\equiv 0 \, (\mod p)$. Finally, if $p^k$ divides $m$ and $n$, then $m=p^km'$, $n=p^kn'$ and
\begin{align*}
\pi_k\Big( X_i\at{m}X_i\at{n} &- \binom{m+n}{m}X_i\at{m+n} \Big)= X_i\at{m'}X_i\at{n'} - \binom{p^k(m'+n')}{p^km'}X_i\at{m'+n'} \\
&= X_i\at{m'}X_i\at{n'} - \binom{m'+n'}{m'}X_i\at{m'+n'}=0.
\end{align*}
The proof for the $Y_i\at{m}$'s is analogous.

Notice that $\pi_k(X_i\at{m}Y_j\at{n}-Y_j\at{n}X_i\at{m})=0$ if $i\neq j$. For \eqref{eq:reln Xn Ym}, it is enough
to verify that \eqref{eq:reln Xn Ym, case} is annihilated since $X_i\at{M}$, $Y_i\at{N}$ can be written as products
of $X_i\at{p^m}$, $Y_j\at{p^n}$. If either $m<k$ or else $n<k$, then $\pi_k$ annihilates both sides of the equality.
Let $m=n=k$. Then
\begin{align*}
 \pi_k\Big( [X_i\at{p^k},Y_i\at{p^k}] &-\sum_{\ell=1}^{p^k}  Y_i\at{p^k-\ell} \left( \sum_{t=0}^\ell \binom{\ell+t}{\ell-t} H_i\at{t} \right) X_i\at{p^k-\ell}\Big) \\
 &= [X_i,Y_i]-H_i=0.
\end{align*}
Now $\pi_k$ annihilates both equations of \eqref{eq:reln Hn Xm} by direct computation.

For \eqref{eq:rels dist-alg -9}, $\pi_k$ annihilates the left hand side if $p^k$ does not divide either $m$ or else $n$.
If $m=p^km'$, $n=p^kn'$ with $n>-ma_{ij}$, then $n'>-m'a_{ij}$ and
\begin{align*}
 \pi_k \left( \ad \left( X_i\at{n}\right)\left( X_j\at{m}\right)\right) = \ad \left( X_i\at{n'}\right)\left( X_j\at{m'}\right)=0.
\end{align*}
Finally \eqref{eq:rels dist-alg -10} follows analogously. Hence $\pi_k$ is an algebra map.

To see that $\pi_k$ is a Hopf algebra map, it remains to prove that $\pi_k$ is a coalgebra map. But it follows
since the elements $X_i\at{p^j}$, $Y_i\at{p^j}$, $H_i\at{p^j}$, $0\le j\le k$, which are primitive elements and generate $\cD_{k+1}{G}$
as an algebra, are applied to primitive elements of $\cD_{1}{G}$.
\epf

The map $\pi_k$ fits in an exact sequence of Hopf algebras.

\begin{prop}\label{prop:extn Hopf algs}
The sequence of Hopf algebras
\begin{align}\label{eq:ex seq Dist}
\xymatrix{ \Bbbk \ar@{->}[r] & \cD_{k}{G} \ar@{^{(}->}[r] & \cD_{k+1}{G} \ar@{->>}[r]^{\pi_k} & \cD_{1}{G} \ar@{->}[r] & \Bbbk }
\end{align}
is exact.
\end{prop}
\pf
By Lemmas \eqref{lemma:normal Hopf subalg} and \eqref{lemma:surj Hopf alg map}, it remains to prove that
\begin{itemize}
 \item $\ker \pi_k= \cD_{k+1}{G}(\cD_{k}{G})^+$, and
 \item $\cD_{k}{G}=\cD_{k+1}{G}^{\co \pi_k}=\{x\in\cD_{k+1}{G}: (\id\otimes\pi_k)\Delta(x)=x\otimes 1 \}$.
\end{itemize}
Note that $\cD_{k+1}{G}(\cD_{k}{G})^+ \subseteq \ker\pi_k$ since $(\cD_{k}{G})^+$ is spanned by $X_i\at{k}$, $Y_i\at{k}$,
$H_i\at{k}$, $1\le k\le p^{n-1}$; the equality follows because both subspaces have the
same dimension, $\dim\cD_{k+1}{G}-\dim\cD_{k}{G}$. Now $\cD_{k+1}{G}^{\co \pi_k} \supseteq \cD_{k}{G}$, and the
equality follows by \cite[Theorem 3.4]{Tak - HopfGalois}.
\epf

\subsection{Steinberg decomposition for simple modules}

The purpose of this section is to introduce the Steinberg's tensor product Theorem. We will prove an analogous result for our quantum version of algebra of distributions of $SL_2$. In order to state this result, we fix some notation \cite{Jan - book}.

\medspace
Let $T\leq G$ be a maximal split torus and $X = X(T)$ be the group of characters of $T$. 
$R$ is the associated root system, $S$ is a fixed basis of $R$ and $R^+$ is the set of positive roots corresponding to $S$.
For each $\alpha\in R$, let $\alpha^{\vee}$ be the associated coroot, and 
$\langle\beta,\alpha^{\vee}\rangle$ denotes the natural pairing, with the normalization $\langle \alpha, \alpha^{\vee} \rangle = 2$ for all $\alpha\in S$.

We consider the following subsets of $X$:
\begin{align*}
X_+ &= \{\lambda\in X\mid \langle \lambda,\alpha^{\vee} \rangle \geq 0, \ \forall \alpha\in R^+\}, 
\mbox{ the set of dominant weights,} \\
X_r &= \{\lambda\in X_+\mid \langle\lambda,\alpha^{\vee}\rangle < p^r, \ \forall \alpha\in S\},
\text{ the set of $r$-restricted weights}, \ r\geq 1.
\end{align*}
Recall that the assignment $\lambda\mapsto L(\lambda)$ establishes a bijection between $X^+$ and the simple $G$-modules up to isomorphism.

Let $B\leq G$ be the Borel subgroup containing $T$ corresponding to $R^-=-R^+$. Given $\lambda\in X_+$ and $M$ a (rational) $G$-module, we set 
\begin{itemize}[leftmargin=*]
\item $M_{\lambda} = \{m\in M\mid t.m = \lambda(t)m\mbox{ for all }t\in T\}$ is the $\lambda$-weight space;
\item $\nabla(\lambda) = \ind_B^G(\lambda)$, the costandard module of highest weight $\lambda$;
\item $L(\lambda) = \soc_G\nabla(\lambda)$, the simple module with highest weight $\lambda$.
\end{itemize}
Let $\cF: G\to G$ be the Frobenius morphism: it arises from the map $\Bbbk\to \Bbbk$, $x\mapsto x^p$. Then $M^{[r]}$ is the $G$-module over the underlying additive group $M$ with $G$-action obtained up to compose the original $G$-action with $\cF^r$.

\begin{theorem}\label{thm:Steinberg}
\cite[Proposition II.3.16]{Jan - book}
Let $r\in\N$, $\lambda\in X_r$, $\mu\in X_+$. Then
\begin{align*}
L(\lambda+p^r \mu)& \simeq L(\lambda)\otimes L(\mu)^{[r]}.
\end{align*}
\end{theorem}

\section{Some cleft extensions of $\u_\lambda(\Sl_2)$}

We introduce the algebras $\Dl{N}{\Sl_2}$, $N\in\N_0$, and prove some properties about their algebra structure which mimic \S \ref{section:dist alg}.

\subsection{$q$-numbers}

We use the following $q$-numbers as in \cite{L - book},
\begin{align}\label{eq:def qbinomial}
 [m]_\lambda&:= \frac{\lambda^m-\lambda^{-m}}{\lambda-\lambda^{-1}},   &
 \fac{m} &= (m)_\lambda (m-1)_\lambda \dots (1)_\lambda, \\
 \bic{m}{n} &:= \prod_{j=1}^n \frac{\lambda^{m-j+1}-\lambda^{-m+j-1}}{\lambda^j-\lambda^{-j}}, &
 0 & \leq n <\ell.
\end{align}

Let $\lambda$ a primitive root of unity of order $\ell$; we assume that $\ell>1$ is odd.

Now we need $q$-\emph{binomial numbers} associated to the $\ell$-expansion. Set
\begin{align}\label{eq:def gen qbinomial}
 \bil{m}{n} &:= \prod_{i\geq 0} \bic{m_i}{n_i}, &  m &=\sum_{i\geq0}m_i \ell^i, \, n=\sum_{i\geq0}n_i \ell^i, \, 0\leq m_i,n_i<\ell.
\end{align}

\begin{lemma}
 Let $m,n,p\geq 0$. Then
\begin{align} \label{eq:prop qbil-1}
 \bil{m+n}{m}&=\bil{m+n}{n}, \\ \label{eq:prop qbil-2}
 \bil{m+n}{n}\bil{m+n+p}{p}&=\bil{n+p}{n}\bil{m+n+p}{m}.
\end{align}
\end{lemma}
\pf
For \eqref{eq:prop qbil-1}, if $m_i+n_i<\ell$ for all $i$, then $(m+n)_i=m_i+n_i$ and
$$ \bil{m+n}{m} = \prod_{i\geq 0} \bic{m_i+n_i}{n_i} = \bil{m+n}{n}. $$
Otherwise there exists $i\geq 0$ such that $m_i+n_i\geq\ell$, we assume $i$ is minimal with this property.
Thus $(m+n)_i=m_i+n_i-\ell< m_i,n_i$, and both sides are 0.

For \eqref{eq:prop qbil-2}, if $m_i+n_i+p_i<\ell$ for all $i$, then $(m+n+p)_i=m_i+n_i+p_i$ and
$$ \bil{m+n}{n}\bil{m+n+p}{p} = \prod_{i\geq 0} \frac{\fac{m_i+n_i+p_i}}{\fac{m_i}\fac{n_i}\fac{p_i}} = \bil{n+p}{n}\bil{m+n+p}{m}. $$
Otherwise there exists $i\geq 0$ such that $m_i+n_i+p_i\geq\ell$, we assume $i$ is minimal with this property.
\begin{itemize}
 \item If $m_i+n_i, n_i+p_i\geq \ell$, then $\bil{m+n}{n}=\bil{n+p}{n}=0$ since $(m+n)_i=m_i+n_i-\ell<n_i$, $(n+p)_i=n_i+p_i-\ell<n_i$.
 \item If $m_i+n_i\geq \ell>n_i+p_i$, then $\bil{m+n}{n}=\bil{m+n+p}{m}=0$ since $(m+n)_i=m_i+n_i-\ell<n_i$, $(n+p)_i=n_i+p_i$, $(m+n+p)_i=m_i+n_i+p_i-\ell<m_i$.
 \item Finally, if $m_i+n_i, n_i+p_i<\ell$, then $\bil{m+n+p}{p}=\bil{m+n+p}{m}=0$.
\end{itemize}
In all the cases, both sides of \eqref{eq:prop qbil-2} are 0.
\epf

\subsection{The Hopf algebra $\u_\lambda(\Sl_2)$}

Throughout this work $\u_\lambda(\Sl_2)$ denotes the algebra presented by generators $E$, $K$, $F$, and relations
\begin{align}\label{eq:defn-rels-uqsl2-1}
& K^\ell=1, & & KF=\lambda^{-2} \, FK, & &   \\ \label{eq:defn-rels-uqsl2-2}
&E^\ell=F^\ell=0, & & KE=\lambda^2 \, EK, & & EF-FE = \frac{K-K^{-1}}{\lambda-\lambda^{-1}}.
\end{align}
It is slightly different from the small quantum group appearing in \cite{L - mod rep}. Then $\u_\lambda(\Sl_2)$ is a Hopf algebra with coproduct:
\begin{align*}
\Delta(K) &= K\ot K, & \Delta(E)&=E\ot 1+K\ot E, & 
\Delta(F)=F\ot K^{-1} + 1\ot F.
\end{align*}

\medbreak

Let $\u_\lambda^+(\Sl_2)$, respectively $\u_\lambda^0(\Sl_2)$, $\u_\lambda^-(\Sl_2)$, be the subalgebra spanned
by $E$, respectively $K$, $F$. Then $\u_\lambda^0(\Sl_2) \simeq \Bbbk \Z/\ell\Z$ while
$\u_\lambda^{\pm}(\Sl_2)$ are isomorphic to $\Bbbk[x]/<x^\ell>$. The multiplication induces a linear isomorphism
$$ \u_\lambda^-(\Sl_2) \otimes \u_\lambda^0(\Sl_2) \otimes \u_\lambda^+(\Sl_2) \simeq \u_\lambda(\Sl_2) . $$
Thus $\{F^a K^b E^c | 0\le a,b,c<\ell \}$ is a basis of $\u_\lambda(\Sl_2)$ and $\dim \u_\lambda(\Sl_2)=\ell^3$.
To simplify the notation in forthcoming computations, let
\begin{align*}
E\at{a} &= \frac{E^a}{\fac{a}}, & F\at{a}&= \frac{F^a}{\fac{a}}, &
\bic{K;s}{a}&= \prod_{j=1}^a \frac{\lambda^{s-j+1}K-\lambda^{-s+j-1}K^{-1} }{\lambda^j- \lambda^{-j}}.
\end{align*}
By direct computation,
\begin{align}\label{eq:conmut-En-Fn}
E\at{m}F\at{n} &= \sum_{i=0}^{\min\{ m,n\} } F\at{n-i} \bic{K;2i-m-n}{i} E\at{m-i},
&  0\le & m,n <\ell.
\end{align}

Let $\u_\lambda^{\ge 0}(\Sl_2)$ be the subalgebra spanned by $E$ and $K$. For each $0\le z<\ell$ the 1-dimensional representation
$\Bbbk_z$ of $\u_\lambda^0(\Sl_2) \simeq \Bbbk \Z/\ell\Z$ given by $K\mapsto \lambda^z$ can be extended to $\u_\lambda^{\ge 0}(\Sl_2)$ by
$E\mapsto 0$. Let $\cM(z)= \u_\lambda(\Sl_2)\otimes_{\u_\lambda^{\ge 0}(\Sl_2)} \Bbbk_z$: it is a $\u_\lambda(\Sl_2)$-module with basis
$v_j:= F\at{j}\otimes 1$, $0\le j<\ell$, such that for all $0\le m,n<\ell$
\begin{align}\label{eq:verma-uqsl2-Fm-vn}
F\at{m} \cdot v_n &= \bic{m+n}{m} v_{m+n}, & K\cdot v_n &=\lambda^{z-2n} v_n,
\\ \label{eq:verma-uqsl2-Em-vn}
E\at{m} \cdot v_n &= \bic{z+m-n}{m} v_{n-m}.
\end{align}
Here $v_n=0$ if either $n<0$ or $n\ge \ell$. Each module $\cM(z)$ has a maximal proper submodule $\cN(z)$. The quotient $\cL(z)=\cM(z)/\cN(z)$ is simple, and has dimension $z+1$: indeed $(v_i)_{0\le i\le z}$ is a basis of $\cL(z)$. Moreover, the family $\{\cL(z)\}_{0\le z<\ell}$ is a set of representatives of the classes of simple modules up to isomorphism.

\subsection{The cleft extensions $\Dl{N}{\Sl_2}$}

We mimic the definition by generators and relations of the algebra of distributions, but in a \emph{quantized} context.

\begin{definition}\label{defn:Dlambda}
Let $\Dl{N}{\Sl_2}$ be the algebra defined by generators $E\atb{i}$, $F\atb{i}$, $K\atb{i}$, $0\le i\le N$ and relations
\begin{align}
K\atb{i} K\atb{j} &= K\atb{j}K\atb{i}, & \left(K\atb{i}\right)^{\ell}&=1; \label{eq:rel D - comm K} \\
K\atb{i} E\atb{j} &=\lambda^{2\delta_{ij}} E\atb{j}K\atb{i}, &  K\atb{i} F\atb{j} &=\lambda^{-2\delta_{ij}} F\atb{j}K\atb{i};
\label{eq:rel D - act K over E,F}\\
E\atb{i} E\atb{j} &= E\atb{j}E\atb{i}, &  F\atb{i} F\atb{j} &= F\atb{j}F\atb{i}; \label{eq:rel D - comm E - comm F}\\
\left(E\atb{i}\right)^{\ell} &= \left(F\atb{i}\right)^{\ell}=0; &
E\atb{i}F\atb{j}&= F\atb{j}E\atb{i}, \quad j\neq i; \label{eq:rel D - powers E,F}
  \end{align}
\begin{equation}\label{eq:rel D - bracket E,F}
 E\atb{j}F\atb{j} = \sum_{t=0}^{\ell^j} F\at{\ell^j-t}  \bik{K;2t-2\ell^j}{t} E\at{\ell^j-t}.
\end{equation}

Here,
$K\atb{-i}:=(K\atb{i})^{-1}$; for $m=\sum_{i=0}^N m_i\ell^i$, $s=\sum_{i=0}^N s_i\ell^i$, $t=\sum_{i=0}^N t_i\ell^i$, $0\le m_i,s_i, t_i<\ell$,
\begin{align*}
E\at{m} &:= \prod_{i=0}^N  \frac{\left(E\atb{i}\right)^{m_i}}{\fac{m_i}},
&
\bik{K;s}{t} &=\prod_{i=0}^N \bic{K\atb{i};s_i}{t_i},
&
F\at{m} &:= \prod_{i=0}^N  \frac{\left(F\atb{i}\right)^{m_i}}{\fac{m_i}}.
\end{align*}
\end{definition}

\medbreak

\begin{remark}
 $\Dl{N}{\Sl_2}$ is $\Z$-graded, with
\begin{align*}
\deg E\atb{i}&=-\deg F\atb{i}=\ell^i, & \deg K\atb{i}&=0, &  0\le & i\le N.
\end{align*}
\end{remark}

\begin{lemma}\label{lemma:proj over uqsl2}
For each pair $M<N$, there exists a surjective algebra map $\pi_{M,N}:\Dl{N}{\Sl_2} \to \Dl{M}{\Sl_2}$ such that
\begin{align*}
\pi_N(X\atb{i})&=
\begin{cases} X\atb{i-N+M}, & i\ge N-M, \\ 0 & i<N-M, \end{cases}
&  X\in\{E,F,K\}.
\end{align*}
In particular, there exists a surjective algebra map $\pi_N:\Dl{N}{\Sl_2} \to \u_\lambda(\Sl_2)$,
 \begin{align*}
  \pi_N(E\atb{i})&=\delta_{iN}E, & \pi_N(F\atb{i})&=\delta_{iN}F, & \pi_N(K\atb{i})&=K^{\delta_{iN}}, & 0\le & i\le N.
 \end{align*}
\end{lemma}
\pf
Straightforward.
\epf

Let $\Dp{N}{\Sl_2}$, resp. $\Dm{N}{\Sl_2}$, $\Dz{N}{\Sl_2}$, be the subalgebras generated by $E\atb{i}$, resp. $F\atb{i}$, $K\atb{i}$, $0\le i\le N$.
Let $\Dge{N}{\Sl_2}$, resp. $\Dle{N}{\Sl_2}$, be the subalgebras generated by $E\atb{i}$ and $K\atb{i}$, resp. $F\atb{i}$ and $K\atb{i}$.

\begin{remark}
\begin{enumerate}[leftmargin=*,label=\rm{(\alph*)}]
  \item  There exists an algebra antiautomorphism $\phi_N$ of $\Dp{N}{\Sl_2}$ such that $\phi_N(E\atb{i})=F\atb{i}$, $\phi_N(F\atb{i})=E\atb{i}$,
 $\phi_N(K\atb{i})=K\atb{i}$, $0\le i\le N$.
  \item There exists an algebra map $\iota_N:\Dl{N-1}{\Sl_2} \to \Dl{N}{\Sl_2}$ which identifies the corresponding generators. Clearly, $\phi_N \circ \iota_N=\iota_N \circ \phi_{N-1}$.
\end{enumerate}
\end{remark}

\begin{lemma}\label{lemma:Verma-modules}
Let $z=\sum_{i=0}^N z_i\ell^i$, $0\le z_i<\ell$. There exists a $\Dl{N}{\Sl_2}$-module $M(z)$ with basis $(v_t)_{0\le t\le \ell^{N+1}-1}$ such that
\begin{align}\label{eq:Verma-modules-action-1}
E\atb{i}\cdot v_t&=\left\{ \begin{array}{ll} \bra{z_i+1-t_i} v_{t-\ell^i}, & t_i>0, \\ 0, & t_i=0; \end{array} \right.  &
K\atb{i}\cdot v_t&=\lambda^{z_i-2t_i} v_t, \\ \label{eq:Verma-modules-action-2}
F\atb{i}\cdot v_t&= \bra{t_i+1} v_{t+\ell^i}, & 0\le & i\le N.
\end{align}
\end{lemma}
\pf
We simply check that $E\atb{i}, F\atb{i}, K\atb{i}\in\End M(z)$, $0\le i\le N$,
satisfy relations \eqref{eq:rel D - comm K}--\eqref{eq:rel D - bracket E,F}. First equation of \eqref{eq:rel D - comm K} holds since
$K\atb{i} K\atb{j}\cdot v_t=\lambda^{z_i+z_j-2t_i-2t_j} v_t$, while the second follows since $\lambda^\ell=1$.
For the first relation in \eqref{eq:rel D - act K over E,F}, both sides annihilate $v_t$ if $t_j=0$; for $t_j\neq 0$, $(t-\ell^j)_i=t_i-\delta_{ij}$, so
\begin{align*}
K\atb{i} E\atb{j}\cdot v_t &= \bra{z_j+1-t_j} \lambda^{z_i-2(t-\ell^j)_i} v_{t-\ell^j}= \lambda^{2\delta_{ij}}
E\atb{j}K\atb{i} \cdot v_t.
\end{align*}
For the second relation, both sides annihilate $v_t$ if $t_j=\ell-1$; for $t_j<\ell-1$,
\begin{align*}
K\atb{i} F\atb{j}\cdot v_t &= \bra{t_j+1} \lambda^{z_i-2(t+\ell^j)_i} v_{t+\ell^j}= \lambda^{-2\delta_{ij}}
F\atb{j}K\atb{i} \cdot v_t,
\end{align*}
since $(t+\ell^j)_i=t_i+\delta_{ij}$. For the first equation in \eqref{eq:rel D - comm E - comm F}, if $t_it_j\neq 0$, $i\neq j$, then
\begin{align*}
E\atb{i} E\atb{j}\cdot v_t &= \bra{z_j-t_j+1} \bra{z_i-(t-\ell^j)_i+1} v_{t-\ell^i-\ell^j} \\
& = \bra{z_j-t_j+1} \bra{z_i-t_i+1} v_{t-\ell^i-\ell^j} = E\atb{j} E\atb{i} \cdot v_t,
\end{align*}
while for $t_it_j=0$, both sides are $0$. The second equation follows similarly.

For the first part of \eqref{eq:rel D - powers E,F}, $(E\atb{i})^{t_i+1}\cdot v_t=(F\atb{i})^{\ell-t_i}\cdot v_t=0$, so
$(E\atb{i})^{\ell}$, $(F\atb{i})^{\ell}$ are 0 as operators on $M(z)$. For the second equality, fix $i\neq j$. If $t_i=0$,
then either $(t+\ell^j)_i=t_i$ or else $t_j=\ell-1$; in any case, $E\atb{i}F\atb{j}\cdot v_t=0=F\atb{j}E\atb{i}\cdot v_t$.
If $t_j=\ell-1$, then again both sides are $0$. Finally set $t_i\neq 0$, $t_j\neq \ell-1$. Hence,
\begin{align*}
E\atb{i}F\atb{j}\cdot v_t &= \bra{t_j+1} \bra{z_i-(t+\ell^j)_i+1} v_{t-\ell^i+\ell^j} \\
& = \bra{(t-\ell^i)_j+1} \bra{z_i-t_i+1} v_{t-\ell^i+\ell^j} = F\atb{j}E\atb{i} \cdot v_t,
\end{align*}

It remains to consider \eqref{eq:rel D - bracket E,F}, which can be written as
\begin{equation}\label{eq:rel D - bracket E,F-v2}
E\atb{j}F\atb{j}-F\atb{j}E\atb{j} -\frac{K\atb{j}-K\atb{-j}}{\lambda-\lambda^{-1}} = \sum_{s=1}^{\ell^j-1} F\at{\ell^j-s}  \bik{K;2s}{s} E\at{\ell^j-s}.
\end{equation}
If $1\le s\le \ell^j-1$, then there exists $i<j$ such that $s_i\neq 0$. If $t_i\geq \ell-s_i$, then
\begin{align*}
\bic{K\atb{i};2s_i}{s_i} (E\atb{i})^{\ell-s_i}\cdot v_t &= \prod_{k=1}^{\ell-s_i} \bra{z_i-t_i+k} \prod_{k=1}^{s_i} \bra{z_i-t_i+k} v_{t-(\ell-s_i)\ell^i}=0.
\end{align*}
If $t_i<\ell-s_i$, then $(E\atb{i})^{\ell-s_i}\cdot v_t=0$. In any case, $\bik{K;2s}{s} E\at{\ell^j-s} \cdot v_t=0$, so the right-hand side of \eqref{eq:rel D - bracket E,F-v2} acts by 0 on each $v_t$. For the left-hand side,
\begin{multline*}
\left( E\atb{j}F\atb{j}-F\atb{j}E\atb{j} -\frac{K\atb{j}-K\atb{-j}}{\lambda-\lambda^{-1}}\right) \cdot v_t \\
= \left( \bra{t_j+1}\bra{z_j-t_j} - \bra{z_j-t_j+1}\bra{t_j}-\bra{z_j-2t_j} \right) v_t=0,
\end{multline*}
when $t_j\neq 0,\ell-1$. If $t_j=0$, then
\begin{multline*}
\left( E\atb{j}F\atb{j}-F\atb{j}E\atb{j} -\frac{K\atb{j}-K\atb{-j}}{\lambda-\lambda^{-1}}\right) \cdot v_t = E\atb{j}\cdot v_{t+\ell^j}-0-\bra{z_j} v_t=0,
\end{multline*}
and finally if $t_j=\ell-1$, then
\begin{align*}
\left( E\atb{j}F\atb{j}-F\atb{j}E\atb{j} -\frac{K\atb{j}-K\atb{-j}}{\lambda-\lambda^{-1}}\right) \cdot v_t
=-\bra{z_j+2}\left( F\atb{j}\cdot v_{t-\ell^j}+ v_t\right)=0.
\end{align*}
In any case, the left-hand side of \eqref{eq:rel D - bracket E,F-v2} also acts by 0 on each $v_t$.
\epf

\begin{lemma}\label{lemma:comodule alg}
There exists an algebra map $\rho_N: \Dl{N}{\Sl_2} \to \u_\lambda(\Sl_2)\otimes \Dl{N}{\Sl_2}$,
\begin{align*}
\rho_N(E\atb{i}) &= 1\otimes E\atb{i}, \quad i<N, & \rho_N(E\atb{N}) &= E \otimes 1 + K\otimes E\atb{N}, \\
\rho_N(F\atb{i}) &= 1\otimes F\atb{i}, \quad i<N, & \rho_N(F\atb{N}) &= F \otimes K\atb{-N} + 1\otimes F\atb{N}, \\
\rho_N(K\atb{i}) &= 1\otimes K\atb{i}, \quad i<N, & \rho_N(K\atb{N}) &= K\otimes K\atb{N}.
\end{align*}
Moreover $\Dl{N}{\Sl_2}$ is a left $\u_\lambda(\Sl_2)$-comodule algebra with this map.
\end{lemma}
\pf
Let $\F$ be the free algebra generated by $E\atb{i}$, $F\atb{i}$, $K\atb{i}$, and $\widetilde\rho_N:\F\to \u_\lambda(\Sl_2)\otimes \Dl{N}{\Sl_2}$
the map defined on the generators as $\rho_N$. We check that $\widetilde\rho_N$ annihilates each defining relation of $\Dl{N}{\Sl_2}$ so it induces the algebra map $\rho_N$. For each relation $\mathtt{r}$ involving only generators $E\atb{i}$, $F\atb{i}$, $K\atb{i}$, $0\le i< N$ we have that $\widetilde\rho_N(\mathtt{r})=1\otimes \mathtt{r}=0$, so we consider those relations involving at least one of the generators $E\atb{N}$, $F\atb{N}$, $K\atb{N}$.

For \eqref{eq:rel D - comm K}, $\widetilde\rho_N \left( (K\atb{N}) ^{\ell} \right)= K^{\ell}\otimes (K\atb{N})^{\ell} = 1 \otimes 1$ and for $i<N$,
$$ \widetilde\rho_N \left( K\atb{i} K\atb{N} -K\atb{N}K\atb{i} \right)= K\otimes \left( K\atb{i} K\atb{N} -K\atb{N}K\atb{i} \right)=0. $$
For \eqref{eq:rel D - act K over E,F} and \eqref{eq:rel D - comm E - comm F}, if $i<N$, then
\begin{align*}
 \widetilde\rho_N (K\atb{i} E\atb{N} &- E\atb{N}K\atb{i})= K\otimes (K\atb{i} E\atb{N} - E\atb{N}K\atb{i})=0, \\
 \widetilde\rho_N (K\atb{N} E\atb{i} &- E\atb{i}K\atb{N})= K\otimes (K\atb{N} E\atb{i} - E\atb{i}K\atb{N})=0, \\
 \widetilde\rho_N (K\atb{N} E\atb{N} &-\lambda^2 E\atb{N}K\atb{N})= (KE-\lambda^2 EK)\otimes K\atb{N} \\
& \qquad + K^2\otimes (K\atb{N} E\atb{N} -\lambda^2 E\atb{N}K\atb{N})=0, \\
\widetilde\rho_N(E\atb{i} E\atb{N} &- E\atb{N}E\atb{i})= K\otimes (E\atb{i} E\atb{N} - E\atb{N}E\atb{i})=0.
\end{align*}
The formulas with $F$ in place of $E$ follow analogously.
For \eqref{eq:rel D - powers E,F},
\begin{align*}
 \widetilde\rho_N \left( (E\atb{N})^{\ell}\right) &= \sum_{j=0}^\ell \bil{\ell}{j} E^{\ell-j} K^j\otimes (E\atb{N})^{j}=0,
\end{align*}
and analogously $\widetilde\rho_N \left( (F\atb{N})^{\ell}\right) =0$.
Finally, for \eqref{eq:rel D - bracket E,F} set $\mathtt{r}_{N}$ as the difference between the two sides of this equation, see also
\eqref{eq:rel D - bracket E,F-v2}. By direct computation,
\begin{align*}
\widetilde\rho_N ( \mathtt{r}_{N} ) & = K \otimes \mathtt{r}_{N}+ \Big(EF-FE- \frac{K-K^{-1}}{\lambda-\lambda^{-1}} \Big) \otimes K\atb{-N} =0.
\end{align*}
Then $\rho_N$ is a well defined algebra map, and gives a left $\u_\lambda(\Sl_2)$-coaction.
\epf

\begin{prop}\label{prop:cleft-ext}
 Let $\rho_N$ as above. Then $\iota_N(\Dl{N-1}{\Sl_2})= ^{\co\rho_N} \Dl{N}{\Sl_2}$, and $\iota_N(\Dl{N-1}{\Sl_2}) \subset \Dl{N}{\Sl_2}$ is a
 $\u_\lambda(\Sl_2)$-cleft extension.
\end{prop}
\pf
Let $\gamma:\u_\lambda(\Sl_2)\to \Dl{N}{\Sl_2}$ be the linear map such that
\begin{align}\label{eq:section}
\gamma(F\at{a} K^b E\at{c})&=\frac{(F\atb{N})^a}{\fac{a}} (K\atb{N})^b \frac{(E\atb{N})^c}{\fac{c}}, & 0\le & a,b,c<\ell.
\end{align}
By direct computation,
\begin{align*}
(\id\otimes\gamma)\circ \Delta(F\at{a} K^b E\at{c}) &= \sum_{i,j} F\at{a-i} K^{b+i+j} E\at{c-j}\otimes \frac{(F\atb{N})^i}{\fac{i}} (K\atb{N})^b \frac{(E\atb{N})^j}{\fac{j}} \\
&= \rho\circ\gamma (F^a K^b E^c),
\end{align*}
so $\gamma$ is map of $\u_\lambda(\Sl_2)$-comodules. We claim that $\gamma$ is convolution invertible. By \cite[Lemma 5.2.10]{Mo},
it is enough to restrict $\gamma$ to the coradical of $\u_\lambda(\Sl_2)$, that is, to $\u^0_\lambda(\Sl_2)$.
Now $\kappa:\u^0_\lambda(\Sl_2)\to \Dl{N}{\Sl_2}$, $\kappa(K^b)= (K\atb{N})^{-b}$, $0\le b<\ell$ is the inverse of
$\gamma_{|\u^0_\lambda(\Sl_2)}$ and the claim follows.

Let $B_N:=\{ F\at{m} K\at{n} E\at{p} | 0\le m,n,p < \ell^{N+1} \}$. We claim that $\Dl{N}{\Sl_2}$ is spanned by $B_N$\footnote{In Proposition
\ref{prop:basis, triang decomp} we shall prove that $B_N$ is indeed a basis of $\Dl{N}{\Sl_2}$}.
Let $I$ be the subspace spanned by $B_N$. Note that $I$ is a left ideal, since it is stable by
left multiplication by $F\atb{n}$, $K\atb{n}$ and $E\atb{n}$ by \eqref{eq:rel D - comm K}-\eqref{eq:rel D - bracket E,F}.
Thus $I=\Dl{N}{\Sl_2}$ since $1\in I$, so $\Dl{N}{\Sl_2}$ is spanned by $B_N$. As
$$ F\at{m} K\at{n} E\at{p}=F\at{m'} K\at{n'} E\at{p'} \frac{(F\atb{N})^{m_N}}{\fac{m_N}} (K\atb{N})^{n_N} \frac{(E\atb{p_N})^c}{\fac{p_N}}, $$
where $0\le m'=m-m_N\ell^N, n'=n-n_N\ell^N, p'=p-p_N\ell^N <\ell^N$, and $F\at{m'} K\at{n'} E\at{p'}\in \iota_N(\Dl{N-1}{\Sl_2})$, we have that
$$ \dim \Dl{N}{\Sl_2}\le \dim \iota_N(\Dl{N-1}{\Sl_2}) \ell^3. $$
As we have a cleft extension, $\Dl{N}{\Sl_2}\simeq ^{\co\rho_N} \Dl{N}{\Sl_2}\otimes \u_\lambda(\Sl_2)$; using this fact and that
$\iota_N(\Dl{N-1}{\Sl_2}) \subset ^{\co\rho_N} \Dl{N}{\Sl_2}$ since $\iota_N$ sends each generator of $\Dl{N-1}{\Sl_2}$ to a coinvariant element, we have that
$$ \dim \Dl{N}{\Sl_2} = \dim ^{\co\rho_N} \Dl{N}{\Sl_2}\ell^3 \ge \dim \iota_N(\Dl{N-1}{\Sl_2}) \ell^3. $$
Hence $\dim ^{\co\rho_N} \Dl{N}{\Sl_2} = \dim \iota_N(\Dl{N-1}{\Sl_2})$, which means that these two subalgebras of $\Dl{N}{\Sl_2}$ coincide.
\epf

\begin{prop}\label{prop:basis, triang decomp}
 \begin{enumerate}[leftmargin=*,label=\rm{(\alph*)}]
  \item\label{item:pos-neg-zero-parts} There exist algebra isomorphisms
  \begin{align*}
   \Dz{N}{\Sl_2} & \simeq \Bbbk(\Z_\ell)^{N+1},  & \Dge{N}{\Sl_2} &\simeq \left(\u_\lambda^{\ge 0}(\Sl_2)\right)^{N+1}, \\
   \Dpm{N}{\Sl_2} & \simeq \left(\u_\lambda^{\pm}(\Sl_2)\right)^{N+1},  & \Dle{N}{\Sl_2} &\simeq \left(\u_\lambda^{\le 0}(\Sl_2)\right)^{N+1}.
  \end{align*}
  \item\label{item:PBW-basis} $B_N:=\{ F\at{m} K\at{n} E\at{p} | 0\le m,n,p < \ell^{N+1} \}$ is a basis of $\Dl{N}{\Sl_2}$.
  \item\label{item:triang-decomp} The multiplication induces a linear isomorphism
  $$ \Dm{N}{\Sl_2} \otimes \Dz{N}{\Sl_2} \otimes \Dp{N}{\Sl_2} \simeq \Dl{N}{\Sl_2} .$$
 \end{enumerate}
\end{prop}
\pf
The algebra $\left(\u_\lambda(\Sl_2)\right)^{N+1}$ is generated by $\Et_i$, $\Ft_i$, $\Kt_i$, $0\le i\le N$,
where each 3-uple $\Et_i$, $\Ft_i$, $\Kt_i$ satisfy \eqref{eq:defn-rels-uqsl2-1}, \eqref{eq:defn-rels-uqsl2-2},
and generators with different subindex commute. There are algebra maps $\Phi^{\ddagger}:\left(\u_\lambda^{\ddagger}(\Sl_2)\right)^{N+1}\to
\mathcal{D}_{\lambda,N}^{\ddagger}(\Sl_2)$, $\ddagger\in\{\pm,0,\ge 0,\le 0\}$, where $\Et_i\mapsto E\atb{i}$, $\Ft_i\mapsto F\atb{i}$,
$\Kt_i\mapsto K\atb{i}$, depending on each case.

For $0\le z<\ell^{N+1}$, let $\Psi_z:\Dl{N}{\Sl_2}\to\End M(z)$ be the algebra map of Lemma \ref{lemma:Verma-modules}.
Notice that $\Psi_z\Phi^-$ is injective, and then $\Phi^-$ is so; thus $\Dm{N}{\Sl_2} \simeq \left(\u_\lambda^{-}(\Sl_2)\right)^{N+1}$.
The map $\Phi^{0}:\Bbbk(\Z_\ell)^{N+1} \to \Dz{N}{\Sl_2}$, $\alpha_i\mapsto K_i$ is surjective. The action of $\Bbbk(\Z_\ell)^{N+1}$ over $v_0$
is given by character $K_i\mapsto \lambda^{z_i}$. Thus $\Bbbk(\Z_\ell)^{N+1} \simeq \Dz{N}{\Sl_2}$. From here we derive that $\Phi^{\le 0}$ is also
an isomorphism. The remaining isomorphisms in \ref{item:pos-neg-zero-parts} follow by using the antiautomorphism $\phi$.

For \ref{item:PBW-basis}, we have to prove that $B_N$ is linearly independent since we have proved that $\Dl{N}{\Sl_2}$ is spanned by $B_N$
in the proof of Proposition \ref{prop:cleft-ext}. We invoke Diamond Lemma \cite[Theorem 1.2]{B - diamond}. Indeed, the lexicographical order
for words written with letters $\{F\atb{i}, K\atb{i},E\atb{i}\}_{0\le i\le N}$ such that
$$ F\atb{0}< \dots< F\atb{N} < K\atb{0}< \dots< K\atb{N} <E\atb{0}< \dots< E\atb{N} $$
is \emph{compatible} (in the notation of loc. cit.) with the reduction system. Each element of $B_N$ is \emph{irreducible}, so $B_N$ is contained
in a basis of $\Dl{N}{\Sl_2}$. Thus $B_N$ is a linearly independent set.
Finally \ref{item:triang-decomp} follows \ref{item:pos-neg-zero-parts} and \ref{item:PBW-basis}.
\epf


\begin{definition}
By Proposition \ref{prop:basis, triang decomp} \ref{item:PBW-basis} each $\iota_N$ is injective. Hence we may consider $\Dl{N-1}{\Sl_2}$ as a subalgebra of $\Dl{N}{\Sl_2}$. Moreover we can consider the inclusions $\iota_{M,N}:\Dl{M}{\Sl_2}\to \Dl{N}{\Sl_2}$ for $M\leq N$, where
\begin{align*}
\iota_{N,N}&=\id_{\Dl{N}{\Sl_2}}, & \iota_{M,N}&=\iota_{M}\iota_{M+1}\dots \iota_{N-1}\mbox{ for }M<N.
\end{align*}
Then we define
\begin{align}\label{eq:defn-Dlambda}
\Dli{\Sl_2} := \lim_{\rightarrow} \Dl{N}{\Sl_2}.
\end{align}
\end{definition}

\section{Finite-dimensional irreducible $\Dl{N}{\Sl_2}$-modules}

Next we study simple modules for the algebras $\Dl{N}{\Sl_2}$. We prove that they are highest weight modules as we can expect, and obtain a decomposition related with the inclusion $\iota_{N-1,N}:\Dl{N-1}{\Sl_2}\to \Dl{N}{\Sl_2}$ and the \emph{Frobenius map} $\pi_{N}:\Dl{N}{\Sl_2}\to \u_\lambda(\Sl_2)$. The tensor product decomposition can be seen as an analogous of Steinberg decomposition, c.f. Theorem \ref{thm:Steinberg}.

\subsection{Highest weight modules}

Now we mimic what is done for simple modules of quantum groups, e. g. \cite[\S 6 \& 7]{L - mod rep}. For the sake of completeness we include the proofs.

\medbreak

Let $V$ be a finite dimensional $\Dl{N}{\Sl_2}$-module. As $\Dz{N}{\Sl_2}$ is the group algebra of $\Z_\ell^{N+1}$, $V$ decomposes as the direct sum of eigenspaces: each $K\atb{i}$ acts by a scalar $\lambda^{p_i}$, $0\le p_i<\ell$. Hence we may encode the data saying that $V=\oplus_{0\le p<\ell^{N+1}} V_p$, where
\begin{align}\label{eq:defn-weight-sp}
	V_p&:= \{ v\in V | K\atb{i}\cdot v=\lambda^{p_i} v \mbox{ for all } 0\le i\le N \}, & p&=\sum_{i=0}^{N}p_i\ell^i.
\end{align}

\begin{definition}
We say that $v\in V$ is a \emph{primitive vector} of weight $p$ if $v\in V_p$ and $E\atb{i}\cdot v=0$ for all $0\le i\le N$. $V$ is called a \emph{highest weight module} if it is generated (as $\Dl{N}{\Sl_2}$-module) by a primitive vector $v$, which is called a \emph{highest weight vector}; its weight $p$ is called a \emph{highest weight}.
\end{definition}

\smallbreak

Given $0\leq p<\ell^{N+1}$, let $\Bbbk_p$ be the 1-dimensional representation of $\Dge{N}{\Sl_2}\simeq \left(\u_\lambda^{\ge 0}(\Sl_2)\right)^{N+1}$ such that $K\atb{i}\cdot 1 =\lambda^{p_i}$ and $E\atb{i}\cdot 1=0$. Let
\begin{align*}
\verma{N}{p} &=\Ind{\Dge{N}{\Sl_2}}{\Dl{N}{\Sl_2}} \Bbbk_p \simeq  \Dl{N}{\Sl_2}\otimes_{\Dge{N}{\Sl_2}} \Bbbk_p.
\end{align*}
Notice that $v_0:= 1\otimes 1\in \verma{N}{p}$ is a primitive vector, and moreover $\verma{N}{p}$ is a highest weight module with highest weight $p$.

\begin{remark}\label{rem:Verma-modules}
Let $v_t=F\at{t}v_0\in \verma{N}{p}$. Then $(v_t)_{0\le t<\ell^{N+1}}$ is a basis of $\verma{N}{p}$, and $\verma{N}{p}$
is isomorphic the module $M(p)$ in Lemma \ref{lemma:Verma-modules}.
Moreover the action on the basis $(v_t)_{0\le t<\ell^{N+1}}$ is given by formulas \eqref{eq:Verma-modules-action-1} and \eqref{eq:Verma-modules-action-2}.

Indeed there is a $\Dge{N}{\Sl_2}$-linear map $\Bbbk_p\to M(p)$ such that $1\mapsto v_0$; it induces a $\Dl{N}{\Sl_2}$-linear map $\verma{N}{p}\to M(p)$, which is surjective by direct computation, and both modules have dimension $\ell^{N+1}$.
\end{remark}

\begin{remark}
Let $V$ a highest weight module of weight $p$. Then each proper submodule is contained in $\oplus_{t\neq p} V_p$; hence $V$ has a maximal proper submodule $\widehat{V}$ and $V/\widehat{V}$ is a simple $\Dl{N}{\Sl_2}$-module, and at the same time a highest weight module of highest weight $p$.
\end{remark}

\begin{definition}
	Let $\irr{N}{p}:= \verma{N}{p}/\widehat{\verma{N}{p}}$; that is, the simple highest weight module obtained as a quotient of $\verma{N}{p}$.
\end{definition}

\begin{prop}\label{prop:simple-modules}
\begin{enumerate}[leftmargin=*,label=\rm{(\alph*)}]
\item\label{item:trivial-action-E} Let $0\le p<\ell^{N+1}$. Then
$$ \{v\in\irr{N}{p}|E\atb{i}v=0\mbox{ for all }0\le i\le N\}=\Bbbk v_0 .$$
\item\label{item:classif-simple-modules} There exists a bijection between $\{p|0\le p<\ell^{N+1}\}$ and the finite-dimensional simple modules of $\Dl{N}{\Sl_2}$ given by $p\mapsto \irr{N}{p}$.
\end{enumerate}
\end{prop}
\pf
\ref{item:trivial-action-E} Let $v\in\irr{N}{p}-0$ be such that $E\atb{i}v=0$ for all $0\le i\le N$.
We may assume that $v$ has weight $t$ for some $0\le t<\ell^{N+1}$, since $E\atb{i}$ applies each eigenspace of the $\Dz{N}{\Sl_2}$ to another.
Thus $v=a\, v_n$ for some $a\in\Bbbk^\times$ and some $0\le n<\ell^{N+1}$, since each 1-dimensional summand in the decomposition
$\verma{N}{p}=\oplus_{0\le n<\ell^{N+1}} \Bbbk v_n$ corresponds to a different eigenspace for the action of
$\Dz{N}{\Sl_2}\simeq \Bbbk (\Z_\ell)^{N+1}$. As $\irr{N}{p}$ is simple, $\irr{N}{p}=\Dl{N}{\Sl_2}v$, but
\begin{align*}
\Dl{N}{\Sl_2} v = \Dle{N}{\Sl_2} v = \Dle{N}{\Sl_2} v_n \subseteq \oplus_{n\le m<\ell^{N+1}} \Bbbk v_m.
\end{align*}
Hence $n=p$ and the claim follows.
\medskip

\ref{item:classif-simple-modules} Let $\cL$ be a simple $\Dl{N}{\Sl_2}$-module. As a $\Dz{N}{\Sl_2}$-module, $\cL=\oplus \cL_t$.
We pick $v\in\cL_t-0$. We may assume that $E\atb{i}v=0$ for all $0\le i\le N$. Indeed, if $E\atb{j}v=0$ for $j=0,\dots, i-1$ but
$E\atb{i}v\neq 0$, let $n\geq 0$ be such that $w:=(E\atb{i})^n v\neq 0$, $(E\atb{i})^{n+1}v=0$. Then $n<\ell$ since $(E\atb{i})^\ell=0$,
and $w$ satisfies $E\atb{j}w=0$ for $j=0,\dots, i$ since $E\atb{j}E\atb{i}=E\atb{i}E\atb{j}$.

Now there exists a $\Dge{N}{\Sl_2}$-linear map $\widetilde{\phi}:\Bbbk_t\to\cL$, $1\mapsto v$, which induces
a $\Dl{N}{\Sl_2}$-linear map $\phi:\verma{N}{t}\to\cL$ such that $1\mapsto v$. As $\cL$ is simple, $\Dl{N}{\Sl_2}v=\cL$,
so $\phi$ is surjective. Hence $\ker \phi\neq 0$ is a proper submodule of $\verma{N}{t}$ and $\cL\simeq \verma{N}{t}/\ker\phi$ is simple.
Thus $\cL\simeq \irr{N}{t}$.

By \ref{item:trivial-action-E}, $\irr{N}{p}\not\simeq\irr{N}{t}$ if $p\neq t$, and the claim follows.
\epf

\subsection{A tensor product decomposition}

\begin{prop}\label{prop:simple-module-trivial}
\begin{enumerate}[leftmargin=*,label=\rm{(\alph*)}]
\item\label{item:restriction-simple-pN=0} Let $0\le p<\ell^{N}$. Then
\begin{align}\label{eq:extension-N-1-to-N}
E\atb{N}\cdot v&=F\atb{N}\cdot v=0, & K\atb{N}\cdot v &=v, &
\mbox{for all }v &\in \irr{N}{p}.
\end{align}
Moreover, $\irr{N}{p}\simeq \irr{N-1}{p}$ as $\Dl{N-1}{\Sl_2}$-modules.

\item\label{item:extension-simple-pN=0} Reciprocally $\irr{N-1}{p}$ may be endowed of an $\Dl{N}{\Sl_2}$-action by extending the $\Dl{N-1}{\Sl_2}$-action via \eqref{eq:extension-N-1-to-N}, and $\irr{N-1}{p} \simeq \irr{N}{p}$ as $\Dl{N}{\Sl_2}$-modules.
\end{enumerate}
\end{prop}
\pf
\ref{item:restriction-simple-pN=0} By the first equation of \eqref{eq:Verma-modules-action-1}, $E\atb{i}v_{\ell^N}=0$
for all $0\le i\le N$, so $\Dl{N}{\Sl_2}v_{\ell^N}=\Dle{N}{\Sl_2}v_{\ell^N}=\oplus_{n\ge \ell^N}\Bbbk v_n$ is a proper submodule
of $\verma{N}{p}$. Hence $v_n=0$ in $\irr{N}{p}$ for all $n\ge\ell^N$, and $\irr{N}{p}$ is spanned by (the image of) $(v_m)_{0\le m<\ell^N}$.
Thus \eqref{eq:extension-N-1-to-N} follows by this fact and \eqref{eq:Verma-modules-action-1}-\eqref{eq:Verma-modules-action-2}.

By \eqref{eq:extension-N-1-to-N}, $W\subset \irr{N}{p}$ is a $\Dl{N-1}{\Sl_2}$-submodule if and only if $W$ is a $\Dl{N}{\Sl_2}$-submodule.
Hence $\irr{N}{p}$ is simple as $\Dl{N-1}{\Sl_2}$-module and the last statement follows.

\medbreak
\ref{item:extension-simple-pN=0} We have to check all the defining relations \eqref{eq:rel D - comm K}-\eqref{eq:rel D - bracket E,F}
of $\Dl{N}{\Sl_2}$. Those not involving $E\atb{N}$, $F\atb{N}$, $K\atb{N}$ follow since $\irr{N-1}{p}$ is a $\Dl{N-1}{\Sl_2}$-module,
and relations $E\atb{N}$, $F\atb{N}$, $K\atb{N}$ follow easily except \eqref{eq:rel D - bracket E,F} for $j=N$.
It is equivalent to \eqref{eq:rel D - bracket E,F-v2}, whose left-hand side acts by $0$ on each $v_t$. For the right-hand side, if
$1\le s\le \ell^N-1$, then there exists $i<N$ such that $s_i\neq 0$, and as in the proof of Lemma \ref{lemma:Verma-modules},
$\bik{K;2s}{s} E\at{\ell^j-s} \cdot v_t=0$, so the right-hand side of \eqref{eq:rel D - bracket E,F-v2} acts by 0 on each $v_t$.
Now $\irr{N-1}{p}$ is a highest weight module as $\Dl{N}{\Sl_2}$-module, with highest weight $p$, and simple at the same time, so
$\irr{N-1}{p} \simeq \irr{N}{p}$ as $\Dl{N}{\Sl_2}$-modules.
\epf

\begin{remark}\label{rem:scalar-ext-ul-Dl}
Thanks to the algebra map $\pi_N:\Dl{N}{\Sl_2} \to \u_\lambda(\Sl_2)$, every $\u_\lambda(\Sl_2)$-module is canonically a $\Dl{N}{\Sl_2}$-module.
In particular each simple $\u_\lambda(\Sl_2)$-module $\cL(p)$, $0\le p<\ell$, is a $\Dl{N}{\Sl_2}$-module.
\end{remark}

\begin{lemma}\label{lem:simple-module-trivial}
Let $p=p_N\ell^N$, $0\le p_N<\ell$. Then $\irr{N}{p}\simeq \cL(p_N)$.
\end{lemma}
\pf
As $\pi_N$ is surjective, $W$ is a $\Dl{N}{\Sl_2}$-submodule of $\cL(p_N)$ if and only if $W$ is a $\u_\lambda(\Sl_2)$-submodule.
Thus $\cL(p_N)$ is a simple $\Dl{N}{\Sl_2}$-module. Now
\begin{align*}
E\atb{i}v_0&=0, & K\atb{i}v_0&=\lambda^{p_N\delta_{iN}}v_0, & \mbox{for all }& 0\leq i\le N.
\end{align*}
Hence $v_0\in\cL(p_N)-0$ is a highest weight vector of weight $p=p_N\ell^N$ and the Lemma follows by Proposition \ref{prop:simple-modules}
\epf

\begin{remark}\label{rem:Dl-ul-tensor-modules}
Recall that $\Dl{N}{\Sl_2}$ is an $\u_\lambda(\Sl_2)$-comodule algebra, so the category of $\Dl{N}{\Sl_2}$-modules is a module category over the category of $\u_\lambda(\Sl_2)$-modules: Given a $\u_\lambda(\Sl_2)$-module $\cM$ and a $\Dl{N}{\Sl_2}$-module $\cN$, $\cM\otimes \cN$ is naturally a $\Dl{N}{\Sl_2}$-module via $\rho$.
\end{remark}

Finally we use Remark \ref{rem:Dl-ul-tensor-modules} to describe a \emph{tensor product decomposition} of simple $\Dl{N}{\Sl_2}$-modules.

\begin{theorem}\label{thm:simple-modules-tensor-decomp}
Let $p=p_N\ell^N+\widehat{p}$, where $0\le \widehat{p}<\ell^N$, $0\le p_N<\ell$. Then
\begin{align*}
\irr{N}{p} & \simeq \cL(p_N)\otimes \irr{N}{\widehat{p}} & \mbox{as }\Dl{N}{\Sl_2}-\mbox{modules.}
\end{align*}
\end{theorem}

\pf
Let $v_0'$, $v_0''$ be highest weight vectors of $\cL(p_N)$, $\irr{N}{\widehat{p}}$, respectively.
We denote $L=\cL(p_N)\otimes \irr{N}{\widehat{p}}$. As $\cL(p_N)$ is generated by $\{v_t'|0\le t<\ell\}$ as in \eqref{eq:verma-uqsl2-Fm-vn}, and $\irr{N}{\widehat{p}}$ is generated by $\{v_t''=F\atb{t}v_0''|0\le t<\ell^N\}$, see Proposition \ref{prop:simple-module-trivial},
$L$ is generated by $\{v_t=v_{t_N}'\otimes v_{\widehat{t}}'' | 0\le t=\widehat{t}+t_N\ell^N <\ell^{N+1}\}$.
Given $F\at{m} K\at{n} E\at{p}\in B_N$, $0\le m,n,p < \ell^{N+1}$,
we may write
\begin{align*}
F\at{m} K\at{n} E\at{p}&=F\at{m_N\ell^N} K\at{n_N\ell^N} E\at{p_N\ell^N}F\at{m'} K\at{n'} E\at{p'}, & 0\le &m',n',p'<\ell^N.
\end{align*}
Here, $F\at{m'} K\at{n'} E\at{p'}\in \Dl{N-1}{\Sl_2}= ^{\co\rho_N} \Dl{N}{\Sl_2}$, cf. Proposition \ref{prop:cleft-ext}. Thus
\begin{align*}
F\at{m} K\at{n} E\at{p}(y\otimes z) & = F\at{m_N}K\at{n_N}E\at{p_N}y \otimes F\at{m'} K\at{n'} E\at{p'} z,
\end{align*}
for all $y\in \cL(p_N)$, $z\in\irr{N}{\widehat{p}}$, where we use \eqref{eq:extension-N-1-to-N}. From here,
$v_0=v_{0}'\otimes v_{0}''$ is a primitive vector, and $L$ is a highest weight module of highest weight $p$. Thus
it suffices to prove that $L$ is simple. Let $W$ be a submodule of $L$. In particular, $W$ is a $\Dz{N}{\Sl_2}$-submodule,
so it decomposes as a direct sum of eigenspaces; each $v_t$, $0\le t<\ell^{N+1}$, spans the eigenspace of weight $t$, so
we may assume that $v_t\in W$ for some $t$. Let $t$ be minimal. Hence
\begin{align*}
0&=E\atb{N}v_t=E \, v_{t_N}'\otimes v_{\widehat{t}}'', & 0&=E\atb{j}v_t=v_{t_N}'\otimes E\atb{j}v_{\widehat{t}}'', \, 0\le j<N,
\end{align*}
so $E \, v_{t_N}'=0=E\atb{j}v_{\widehat{t}}''$, $0\le j<N$. From here, $t_N=\widehat{t}=0$, and then $W=L$.
\epf

\begin{remark}
$\Dl{N}{\Sl_2}$ is an \emph{augmented algebra} via the map $\eps: \Dl{N}{\Sl_2}\to \Bbbk$,
\begin{align}\label{eq:DlN-augmented}
\epsilon(E\atb{j})&=\epsilon(F\atb{j})=0, & \epsilon(K\atb{j})&=1, &
\mbox{for all }0 &\le j\le N.
\end{align}
Thence $\Bbbk$ is a $\Dl{N}{\Sl_2}$-module and $\Bbbk\simeq \irr{N}{0}$ via $\epsilon$, so $\irr{N}{p}\simeq \cL(p_N)\otimes \Bbbk$
if $p=p_N\ell^N$, $0\le p_N<\ell$.
\end{remark}

\section*{Acknowledgement} The main part of this paper was written during my visit to the Max-Planck-Institute in Bonn as an Alexander von Humboldt Fellow. I would like to thank the institute for its excellent working environment and support. I would like to thank especially to Geordie Williamson for all the discussions and the guidance.

\end{document}